\documentclass[11pt]{article}

\usepackage{xypic}
\usepackage{amsgen}
\usepackage{amsmath}
\usepackage{amstext}
\usepackage{amsbsy}
\usepackage{amsopn}
\usepackage{amsfonts}
\usepackage{amssymb}
\usepackage{graphicx}
\usepackage{pstricks}
\xyoption{all}

\def\Box{\square}
\def\edge{\relbar\joinrel\relbar}

\def\tra#1{\smash{\mathop{\mid\kern
-1pt\joinrel\relbar\joinrel\relbar}\limits^{*}_{#1}}}
\def\longtra#1{\smash{\mathop{\mid\kern
-1pt\joinrel\relbar\joinrel\relbar\joinrel\relbar}\limits^{*}_{#1}}}
\def\vlongtra#1{\smash{\mathop{\mid\kern
-1pt\joinrel\relbar\joinrel\relbar\joinrel\relbar\joinrel\relbar}\limits^{*}_{#1}}}
\def\vvlongtra#1{\smash{\mathop{\mid\kern
-1pt\joinrel\relbar\joinrel\relbar\joinrel\relbar\joinrel\relbar\joinrel\relbar}\limits^{*}_{#1}}}
\def\vvvlongtra#1{\smash{\mathop{\mid\kern
-1pt\joinrel\relbar\joinrel\relbar\joinrel\relbar\joinrel\relbar\joinrel\relbar\joinrel\relbar}\limits^{*}_{#1}}}
\def\etra#1{\smash{\mathop{\mid\kern
-1pt\joinrel\relbar\joinrel\relbar}\limits_{#1}}}

\def\iff{\Leftrightarrow}
\def\Rw{\Rightarrow}
\def\oo{\overline}

\def\wt{\widetilde}

\newcommand{\N}{{\rm I}\kern-2pt {\rm N}}

\def\rat{\mbox{Rat}\,}

\def\aut{\mbox{Aut}\,}
\def\endo{\mbox{End}\,}

\def\fix{\mbox{Fix}\,}
\def\per{\mbox{Per}\,}

\def\p{\varphi}

\def\inv{^{-1}}

\def\bi{\begin{itemize}}
\def\ei{\end{itemize}}
\def\beq{\begin{equation}}
\def\eeq{\end{equation}}

\xyoption{all}

\newtheorem{T}{Theorem}[section]
\newcommand{\bt}{\begin{T}}
\newcommand{\et}{\end{T}}
\newcommand{\ftd}{$\square$\end{T}}

\newtheorem{Proposition}[T]{Proposition}
\newcommand{\bp}{\begin{Proposition}}
\newcommand{\ep}{\end{Proposition}}
\newcommand{\fpd}{$\square$\end{Proposition}}

\newtheorem{Lemma}[T]{Lemma}
\newcommand{\bl}{\begin{Lemma}}
\newcommand{\el}{\end{Lemma}}
\newcommand{\fld}{$\square$\end{Lemma}}

\newtheorem{Corol}[T]{Corollary}
\newcommand{\bc}{\begin{Corol}}
\newcommand{\ec}{\end{Corol}}
\newcommand{\fcd}{$\square$\end{Corol}}

\newtheorem{Result}[T]{Result}
\newcommand{\br}{\begin{Result}}
\newcommand{\er}{\end{Result}}
\newcommand{\frd}{$\square$\end{Result}}

\newtheorem{Example}[T]{Example}
\newcommand{\be}{\begin{Example}}
\newcommand{\ee}{\end{Example}}

\newtheorem{Problem}[T]{Problem}
\newcommand{\bq}{\begin{Problem}}
\newcommand{\eq}{\end{Problem}}

\newcommand{\proof}
   {\par\medbreak\noindent{\bf Proof}.\enspace}

\newcommand{\qed}{
$\Box$
\par\bigbreak}

\normalbaselines
\def\abstract#1{\par\bigskip
\begingroup\small
\baselineskip=12truept
\begin{center}ABSTRACT\end{center}
\par\medskip\par\noindent
\null\hfill\hbox{\vbox{\hsize=5truein\noindent#1}}
\hfill\null\par\endgroup\par}

\title{Fixed points of endomorphisms of graph groups}
\author{{\bf Emanuele Rodaro, Pedro V. Silva}\\ $ $\\ {\em Centro de
Matem\'{a}tica, Faculdade de Ci\^{e}ncias, Universidade do
Porto,}\\ {\em R. Campo Alegre 687, 4169-007 Porto, Portugal}\\
{\em e-mail:} emanuele.rodaro@fc.up.pt, pvsilva@fc.up.pt\\
$ $\\
{\bf Mihalis Sykiotis}\\ $ $\\ {\em Department of Mathematics,
  National and Kapodistrian University of Athens,}\\
{\em Panepistimioupolis, GR-157 84, Athens, Greece}\\
{\em e-mail:} msykiot@math.uoa.gr}
\date{\today}

\begin{document}
\maketitle

\begin{center}\small

\medskip

\end{center}

\abstract{It is shown, for a given graph group $G$, that the fixed
  point subgroup Fix$\,\p$ is finitely
  generated for every endomorphism $\p$ of $G$
  if and only if $G$ is a free product of free abelian
  groups. The same conditions hold for the subgroup of periodic
  points. Similar results are obtained for automorphisms, if the
  dependence graph of $G$ is a transitive forest.}

\section{Introduction}

Gersten proved in the 
eighties that the fixed point subgroup of an automorphism of a free
group of finite rank is always
finitely generated \cite{Ger}.
Using a different approach, Cooper gave an alternative topological
proof \cite{Coo}. This result was generalized to further classes of groups and
endomorphisms in subsequent years. Goldstein and Turner extended it to
monomorphisms of free groups \cite{GT}, and later to
arbitrary endomorphisms 
\cite{GT2}. Collins and Turner extended it to automorphisms of free
products of freely indecomposable groups \cite{CT} (see the survey by
Ventura \cite{Ven}). With respect to automorphisms, the widest
generalization is to hyperbolic groups and is due to Paulin \cite{Pau}.

These results inspired Bestvina and Handel to develop a new line of
research through their innovative train track techniques,
bounding the rank of the fixed point subgroup for the free group
automorphism case \cite{BH}. These results were generalized by the
third author by considering symmetric
endomorphisms of free products and using the concept of Kurosh rank
\cite{Syk,Syk2,Syk3}. These results are general enough to imply
finiteness theorems for arbitrary endomorphisms of finitely generated
virtually free groups.

In \cite{Sil2}, in the sequence of previous work developped with Cassaigne
\cite{CS,CS2}, the second author used automata-theoretic methods to
obtain finiteness results for monomorphisms of monoids defined by
special confluent rewriting systems, and in \cite{Sil3},
followed a different approach that 
also implied the virtually free group endomorphism case.

In the present paper, we discuss which graph groups (also widely known
as right angled Artin groups) admit finitely generated fixed point
(and also periodic point) 
subgroups for every endomorphism/automorphism. A complete solution is
reached in the endomorphism case, a partial solution in the
automorphism case. 

The fact that not all graph groups satisfy these properties (Theorems
\ref{fixpoints} and \ref{fixpointsaut}) shows
that the behaviour of fixed point subgroups imposes strong restrictions on the
structure of the ambient group.

\section{Preliminaries}

Given a group $G$, we denote by $\aut G$ (respectively $\endo G$) the
automorphism group (respectively endomorphism monoid) of
$G$. 
Given $\p \in \endo G$, we say that $g \in G$ is a {\em fixed point}
of $\p$ if $g\p = g$. If $g\p^n = g$ for some $n \geq 1$, we say that
$g$ is a {\em periodic point} of $\p$. Let
$\fix\p$ (respectively $\per\p$) denote the set of all fixed points
(respectively periodic points) of $\p$. Clearly,
\beq
\label{perfix}
\per\p = \displaystyle\bigcup_{n \geq 1} \fix\p^n.
\eeq

The {\em free group} on an alphabet $A$ is denoted by $F_A$. We denote
by $R_A$ the set of all reduced words on $\wt{A} = A \cup
A\inv$. For every $g \in F_A$, we denote by $\oo{g}$ the unique
reduced word representing $g$.

We need to introduce {\em rational subsets} for an arbitrary monoid
$M$. Given $X \subseteq M$, let $X^*$ denote the submonoid of $M$
generated by $X$.
We denote by $\rat M$ the smallest class ${\cal{X}}$ of subsets of $M$
containing the finite sets and satisfying
$$X,Y \in {\cal{X}} \; \Rw \; X\cup Y, XY, X^* \in {\cal{X}}.$$
The following classical result, known as Benois Theorem, relates
through reduction the rational subsets of the free group $F_A$ with the
rational subsets of the free monoid $\wt{A}^*$:

\bt
\label{benois}
{\rm \cite{Ben}}
Let $X \subseteq F_A$. Then $X \in$ {\rm
  Rat}$\, F_A$ if and only if $\oo{X} \in$ {\rm
  Rat}$\, \widetilde{A}^*$.
\et

As a consequence, we have:

\bc
\label{bencon}
{\rm \cite{Ben}}
{\rm Rat}$\, F_A$ is closed under
the boolean operations.
\ec

A (finite) {\em independence alphabet} is an ordered pair of the form $(A,I)$,
where $A$ is a (finite) set and $I$ is a symmetric anti-reflexive relation on
$A$. We can view $(A,I)$ as an undirected graph without loops or
multiple edges, denoted by $\Gamma(A,I)$, by taking $A$ as the vertex
set and $I$ as the edge 
set. Conversely, every such graph determines an independence alphabet.

Let $N(A,I)$ denote the
normal subgroup of $F_A$ generated by the commutators
$$\{ [a,b] \mid (a,b) \in I \}.$$
The {\em graph group} $G(A,I)$ is the quotient $F_A/N(A,I)$, i.e the
group defined by the group presentation 
$$\langle A \mid [a,b] \; ((a,b) \in I) \rangle.$$
Such groups are also known as {\em right angled Artin groups}, {\em
  free partially commutative groups} or even {\em trace groups}. 

If $I = \emptyset$, then $G(A,I)$ is the free group on $A$. Let
$\Delta_A = \{ (a,a) \mid a \in A \}$. If $I = (A
\times A) \setminus \Delta_A$ (i.e $\Gamma(A,I)$ is a complete
graph), then $G(A,I)$ is the free abelian group
on $A$. It is well known that the class of (finitely generated) graph
groups is closed under free product and direct product: the free product
$G(A_1,I_1) \ast G(A_2,I_2)$ is realized by the disjoint union of the
independence graphs $(A_1,I_1)$ and $(A_2,I_2)$, and the 
direct product $G(A_1,I_1) \times G(A_2,I_2)$ is realized by the graph
obtained from the disjoint union of $(A_1,I_1)$ and $(A_2,I_2)$ by
connecting every vertex 
in $A_1$ to every vertex in $A_2$. However, not all
graph groups can be built from infinite cyclic groups using these two
operators, the simplest example being given by the graph
$$a \edge b \edge c \edge d$$
We remark also that, for every $a \in A$, we may define a homomorphism
$\pi_a:G(A,I) \to \mathbb{Z}$ by $a\pi_a = 1$ and $b\pi_a = 0$ for $b
\in A \setminus \{ a \}$.

Finally, we need some concepts and results involving free
products. Assume that $G = G_1 \ast G_2 \ast \ldots \ast G_n$. 
By the Kurosh Subgroup Theorem \cite{Kur}, every subgroup of $G$ is a free
product of the form $F_X \ast H_1 \ast \ldots \ast H_m$, where each
$H_j$ is a conjugate of a subgroup of some
$G_i$. The {\em Kurosh rank} of $H$ is the sum $|X| + m$. In 2007, the
third author proved the following result:

\bt
\label{fkr}
{\rm \cite[Theorem 7]{Syk3}}
Let $G = G_1 \ast \ldots \ast G_n$ be a free product of finitely
generated nilpotent and finite groups and let $\p \in$ {\rm End}$\,
G$. Then {\rm Fix}$\,\p$ has Kurosh rank at most $n$.
\et

We say that a group $G$ satisfies the {\em maximal condition on
  subgroups} if every subgroup of $G$ is finitely generated. The
following result, proved by the third author in 2005, will be used in
the discussion of periodic points:

\bt
\label{stat}
{\rm \cite[Corollary 6.3]{Syk2}}
Let $G = G_1 \ast \ldots \ast G_n$ be a free product of groups satisfying
the maximal condition on subgroups. If $H_1 \subseteq H_2 \subseteq
\ldots$ is an ascending chain of subgroups whose
Kurosh ranks are bounded by a natural number $N$, then there exists an
index $m$ such that $H_i = H_m$ for 
every $i \geq m$.
\et

\section{Endomorphisms}

We can reach a complete answer in the case of fixed/periodic points of
endomorphisms:

\bt
\label{fixpoints}
Let $(A,I)$ be a finite independence alphabet. Then the following
conditions are equivalent:
\bi
\item[(i)] {\rm Fix}$\,\p$ is finitely generated for every $\p \in$
  {\rm End}$\, G(A,I)$;
\item[(ii)] {\rm Per}$\,\p$ is finitely generated for every $\p \in$
  {\rm End}$\, G(A,I)$;
\item[(iii)] $I \cup \Delta_A$ is transitive;
\item[(iv)] $\Gamma(A,I)$ is a disjoint union of complete graphs;  
\item[(v)] $G(A,I)$ is a free product of finitely many free abelian
  groups of finite rank.
\ei
\et

\proof
(i) $\Rw$ (iii). Suppose that $I \cup \Delta_A$ is not
transitive. Then there exist $a,b,c \in A$ such that $(a,b),(b,c) \in
I \cup \Delta_A$ and $(a,c) \notin I \cup \Delta_A$. Note that this
implies that $a,b,c$ are all distinct.

Write $G = G(A,I)$.
We define $\p \in \endo G$ by $a\p = ab$, $b\p = b$, $c\p = b\inv c$ and
$x\p = 1$ for $x \in A \setminus \{ a,b,c \}$. We show that
\beq
\label{fixpoints1}
u \in \fix\p \; \Rw \; u\pi_a = u\pi_c.
\eeq
Indeed, we have $u\p\pi_b = u\pi_a + u\pi_b - u\pi_c$, hence $u\p = u$
yields $u\pi_a = u\pi_c$ and (\ref{fixpoints1}) holds.

Write $A_0 = \{ a,c\}$
and let $\Psi:G \to F_{A_0}$ be the homomorphism defined by
$a\Psi = a$, $c\Psi = c$ and $x\Psi = 1$ for $x \in A \setminus \{ a,c
\}$. We claim that 
\beq
\label{fixpoints2}
a^ic^j \in (\fix\p)\Psi \; \iff \; i = j
\eeq
holds for all $i,j \geq 0$.

Indeed, the direct implication follows from
(\ref{fixpoints1}). On the other hand, we have $(a^ic^i)\p =
(ab)^i(b\inv c)^i = a^ic^i$, hence $a^ic^i = 
(a^ic^i)\Psi \in (\fix\p)\Psi$ for every $i \geq 0$. Thus
(\ref{fixpoints2}) holds.

Suppose that $\fix\p$ is finitely generated. Then $(\fix\p)\Psi$ is
a finitely generated subgroup of $F_{A_0}$ and therefore $(\fix\p)\Psi
\in \rat F_{A_0}$. But $a^*c^* \in \rat F_{A_0}$ as well and so by
Corollary \ref{bencon}
also $L = (\fix\p)\Psi \cap a^*c^* \in \rat F_{A_0}$. Now
(\ref{fixpoints2}) yields $L = \{ a^nc^n \mid n \geq 0 \}$. Since $L
\subseteq R_{A_0}$, it follows from Theorem \ref{benois} that $L \in
\rat A_0^*$, a contradiction since this is a famous example of a
non rational language (it is easy to check that it fails the {\em
Pumping Lemma} for rational languages \cite[Lemma
I.4.5]{Ber}). Therefore $\fix\p$ is not finitely generated.

(ii) $\Rw$ (iii).
We adapt the proof of (i) $\Rw$ (iii) taking the same endomorphism
$\p$ and showing that we may replace
$\fix\p$ by $\per\p$ in both (\ref{fixpoints1}) and
(\ref{fixpoints2}). 

Indeed, for every $n \geq 1$, we have 
 $a\p^n = ab^n$, $b\p^n = b$, $c\p^n = b^{-n} c$ and
$x\p^n = 1$ for $x \in A \setminus \{ a,b,c \}$.
Since $u\p^n\pi_b = n(u\pi_a) + u\pi_b - n(u\pi_c)$, then $u\p^n = u$
yields $u\pi_a = u\pi_c$ and (\ref{perfix}) yields
\beq
\label{perpoints1}
u \in \per\p \; \Rw \; u\pi_a = u\pi_c.
\eeq
Now, similarly to the fixed point case, we show that
\beq
\label{perpoints2}
a^ic^j \in (\per\p)\Psi \; \iff \; i = j
\eeq
holds for all $i,j \geq 0$.
Finally, we use (\ref{perpoints1}) and
(\ref{perpoints2}) in similar fashion to reach the desired contradiction.

(iii) $\Rw$ (iv) $\Rw$ (v). Immediate.

(v) $\Rw$ (i). 
Since abelian groups are nilpotent, $G$ is the free product of
finitely many finitely generated nilpotent groups. It follows from Theorem
\ref{fkr} that $\fix\p$ has finite Kurosh rank. Since each of the free
factors in the Kurosh decomposition of $\fix\p$ is itself finitely
generated, it follows that $\fix\p$ is finitely generated.

(v) $\Rw$ (ii). 
Assume that $G(A,I)$ is a free product of $n$ free abelian
  groups of finite rank.
It is well known that every finitely generated abelian group
satisfies the maximal chain condition on subgroups. Consider the chain
of subgroups of $G(A,I)$
\beq
\label{chs}
\fix\p \subseteq \fix\p^{2!} \subseteq \fix\p^{3!} \subseteq
\ldots
\eeq
By Theorem \ref{fkr}, each $\fix\p^{k!}$ has Kurosh rank at most $n$,
hence by Theorem \ref{stat}, there exists some $m$ such that
$\fix\p^{k!} = \fix\p^{m!}$ for every $k \geq m$. By
(\ref{perfix}), and since $\fix\p^r \subseteq \fix\p^s$ whenever
$r|s$, we get
$$\per\p = \displaystyle\bigcup_{k \geq 1} \fix\p^k = \bigcup_{k \geq
  1} \fix\p^{k!} = \fix\p^{m!},$$
therefore $\per\p$ is finitely generated by the implication (v) $\Rw$ (i).
\qed

\section{Automorphisms}

If we restrict our attention to automorphisms, we can prove an
analogue of Theorem \ref{fixpoints}, but only for {\em transitive
  forests}. We say that a graph is a transitive forest if it has no induced
subgraph of one of the following forms:
$$
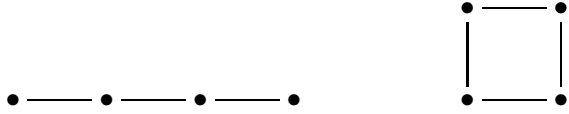
\xymatrix{
&&&&& \bullet \ar@{-}[r] & \bullet \\
\bullet \ar@{-}[r] & \bullet \ar@{-}[r] & \bullet \ar@{-}[r] & \bullet
&& \bullet \ar@{-}[r] \ar@{-}[u] & \bullet \ar@{-}[u]
}$$
Transitive forests constitute an important class of graphs in the
context of graph 
groups since they often establish the territory of positive
algorithmic properties. For instance, results of Lohrey and Steiberg
show that a graph group $G(A,I)$ has solvable
submonoid membership problem (or solvable rational subset membership
problem) if and only if $\Gamma(A,I)$ is a transitive forest \cite{LoS}.

\bt
\label{fixpointsaut}
Let $(A,I)$ be a finite independence alphabet such that $\Gamma(A,I)$
is a transitive 
forest. Then the following conditions are equivalent:
\bi
\item[(i)] {\rm Fix}$\,\p$ is finitely generated for every $\p \in$
  {\rm Aut}$\, G(A,I)$; 
\item[(ii)] {\rm Per}$\,\p$ is finitely generated for every $\p \in$
  {\rm Aut}$\, G(A,I)$;
\item[(iii)] $I \cup \Delta_A$ is transitive;
\item[(iv)] $\Gamma(A,I)$ is a disjoint union of complete graphs;  
\item[(v)] $G(A,I)$ is a free product of free abelian groups.
\ei
\et

\proof
In view of Theorem \ref{fixpoints}, it suffices to show that (i) $\Rw$
(iii) and (ii) $\Rw$ (iii). 

(i) $\Rw$ (iii).
Suppose that $I \cup \Delta_A$ is not
transitive. Then there exist (distinct) $a,b,c \in A$ such that $(a,b),(b,c) \in
I \cup \Delta_A$ and $(a,c) \notin I \cup \Delta_A$. 

Write $G = G(A,I)$.
We define $\p \in \endo G$ by $a\p = ab$, $c\p = b\inv c$ and
$x\p = x$ for $x \in A \setminus \{ a,c \}$. We must show that $\p$ is
well defined, so take $d \in A \setminus \{ a,b,c \}$. All the other
cases being easily checked, it suffices to
show that $(a,d) \in I$ implies $abd = dab$ in $G$ (the
implication $(c,d) \in I \Rw b\inv cd = db\inv c$ in $G$ is
analogous). Indeed, $(a,d) \in I$ implies $(b,d) \in I$ because
$(A,I)$ is a transitive forest: since $(a,c) \notin I$, this is the
only way of avoiding the forbidden configurations. Therefore $\p$ is
well defined.

Now $\p$ is clearly an epimorphism. Green's theorem on graph products
\cite[Corollary 5.4]{Gre} implies that $G$ is residually finite,
and therefore hopfian by a theorem of Mal'cev \cite{Mal}. This implies
that the epimorphism $\p$ is indeed an automorphism. 

Now we use the same argument as in the proof of Theorem
\ref{fixpoints} to show that $\fix \p$ is not finitely generated.

(ii) $\Rw$ (iii). Analogous.
\qed

The obvious question is now what happens in the automorphism case if
$\Gamma(A,I)$ is not a transitive forest. Note that this
excludes $\Gamma(A,I)$ from being a disjoint union of complete graphs!
The following examples show that it can go either way.

\be
\label{fgyes}
Let $\Gamma(A,I)$ be the graph 
$$
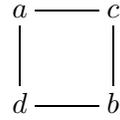
\xymatrix{
a \ar@{-}[r] & c \\
d \ar@{-}[r] \ar@{-}[u] & b \ar@{-}[u]
}$$
Then {\rm Fix}$\,\p$ and {\rm Per}$\,\p$ are both finitely generated
for every $\p \in$ {\rm 
  Aut}$\, G(A,I)$. 
\ee

Indeed, we may write $G = G(A,I) = F_{a,b} \times F_{c,d}$. Given $\p_1
\in \aut F_{a,b}$ and $\p_2
\in \aut F_{c,d}$, let $(\p_1,\p_2) \in \aut G$ be given by
$(u,v)(\p_1,\p_2) = (u\p_1,v\p_2)$ (type I automorphism). Let
$\psi:F_{a,b} \to F_{c,d}$ be 
the isomorphism defined by $a\psi = c$ and $b\psi = d$. Finally, let
$\sigma:G \to F_{c,d} \times F_{a,b}$ be the isomorphism given by
$(u,v)\sigma = (v,u)$. It is immediate that
$(\p_1,\p_2)\sigma(\psi\inv,\psi) \in \aut G$ as well (type II
automorphism). We claim that all the automorphisms of $G$ must be of
type I or type II.

Let $H = (F_{a,b} \times \{ 1 \}) \cup ( \{ 1 \} \times F_{c,d})$. If
$F$ is a free group of rank $> 1$ and $u \in F$, it is well known
\cite{LS} that
the centralizer $C_u$ satisfies
$$C_u \cong \left\{
\begin{array}{ll}
\mathbb{Z}&\mbox{ if }u \neq 1\\
F&\mbox{ if }u = 1
\end{array}
\right.$$
It follows that, for every $(u,v) \in G$,
$$C_{(u,v)} \cong \left\{
\begin{array}{ll}
\mathbb{Z} \times \mathbb{Z} &\mbox{ if }u,v \neq 1\\
G&\mbox{ if }u = v = 1\\
\mathbb{Z} \times F_{a,b}&\mbox{ otherwise}
\end{array}
\right.$$
and so $H$, the set of elements with nonabelian centralizer,
 must be fixed by every $\p \in \aut G$.

Suppose first that $a\p \in F_{a,b} \times \{ 1 \}$. Since $(ab)\p
\neq (ba)\p$, we see that $b\p \in F_{a,b} \times \{ 1 \}$. It 
follows easily that both $F_{a,b} \times \{ 1 \}$ and $\{ 1 \} \times
F_{c,d}$ must be fixed by $\p$. Thus $\p$ must be of type I.

Assume now that $a\p \in \{ 1 \} \times
F_{c,d}$. Similar arguments show that $\p$ admits as restrictions
isomorphisms $\p':F_{a,b} \times \{ 1 \} \to \{ 1 \} \times
F_{c,d}$ and $\p'': \{ 1 \} \times
F_{c,d} \to F_{a,b} \times \{ 1 \}$. Now $\p'(1,\psi\inv)\sigma \in
\aut(F_{a,b} \times \{ 1 \})$, hence $\p'(1,\psi\inv)\sigma =
(\p_1,1)$ for some $\p_1 \in \aut F_{a,b}$. Similarly, $\p''(\psi,1)\sigma =
(1,\p_2)$ for some $\p_2 \in \aut F_{c,d}$. Hence
$$\begin{array}{lll}
(u,v)\p&=&((u,1)\p')((1,v)\p'') =
((u,1)(\p_1,1)\sigma(1,\psi))((1,v)(1,\p_2)\sigma(\psi\inv,1))\\
&=&((u,1)(\p_1,\p_2)\sigma(\psi\inv,\psi))
((1,v)(\p_1,\p_2)\sigma(\psi\inv,\psi))\\
&=&(u,v)(\p_1,\p_2)\sigma(\psi\inv,\psi)
\end{array}$$
holds for every $(u,v) \in G$ and so $\p$ is type II.

We show next that $\fix\p$ is finitely generated for all type I
and type II automorphisms. If $\p = (\p_1,\p_2)$ is type I, then
$\fix\p = \fix\p_1 \times \fix\p_2$ is finitely generated in view of
Gersten's theorem (see also Theorem \ref{fixpointsaut}). Hence we may
assume that $\p = (\p_1,\p_2)\sigma(\psi\inv,\psi)$ is type II. 

Given $(u,v) \in G$, we have $(u,v)\p = (u,v)$ if and only if
$v\p_2\psi\inv = u$ and $u\p_1\psi = v$. This is equivalent to 
$$u\p_1\psi\p_2\psi\inv = u \hspace{.3cm}\wedge\hspace{.3cm}u\p_1\psi
= v.$$
Let $K = \fix \p_1\psi\p_2\psi\inv$. By Gersten's theorem,
$K$ is finitely generated, say $K = \langle z_1,\ldots, z_r
\rangle$. Hence 
$$\fix\p = \{ (u,u\p_1\psi) \mid u \in
K\} = \langle (z_1,z_1\p_1\psi), \ldots, (z_r,z_r\p_1\psi) \rangle$$
is also finitely generated in this case.

Regarding $\per\p$, the type I case is analogous and the type II case
may be reduced to the type case I in view of $\per\p = \per\p^2$.

\be
\label{fgno}
Let $\Gamma(A,I)$ be the graph 
$$\xymatrix{
a \ar@{-}[r] & b \ar@{-}[dr] & \\
d \ar@{-}[r] & c \ar@{-}[u] \ar@{-}[r] & e
}$$
Then there exists some $\p \in$ {\rm 
  Aut}$\, G(A,I)$ such that neither {\rm Fix}$\,\p$ nor {\rm
  Per}$\,\p$ are finitely generated. 
\ee

We adapt the proof of Theorem \ref{fixpoints}.
Define $\p \in \endo G$ by $a\p = ab\inv$, $b\p=b$, $c\p=c$,
$d\p=dc\inv$ and $e\p=ebc$. This is a well-defined homomorphism which
turns out to be an epimorphism, and therefore an automorphism since
$G$ is hopfian. Applying $\pi_b$ and $\pi_c$ to fixpoints yields
\beq
\label{eq:projections_bce}
u \in \fix\p \; \Rw \; u\pi_a = u\pi_d = u\pi_e.
\eeq
Let $A_0 = \{ a,d,e \}$
and let $\Psi:G \to F_{A_0}$ be the homomorphism defined by
$a\Psi = a$, $d\Psi = d$, $e\Psi = e$ and $x\Psi = 1$ for $x \in A
\setminus  A_0$.  Then we use (\ref{eq:projections_bce}) to show that
$$a^ie^jd^k \in (\fix\p)\Psi \; \iff \; i = j = k$$
holds for all $i,j,k \geq 0$.
Intersecting $\fix\p$ with $a^*e^*d^*$ yields a non rational
language (not even {\em context-free}, actually -- see \cite{Ber}) and
a straightforward 
adaptation of the argument used previously proves that
$\fix\p$ is not finitely generated. The discussion of $\per\p$ is also
adapted as in the proof of Theorem \ref{fixpoints}.

\section*{Acknowledgements}

The first two authors acknowledge support from the European Regional Development
Fund through the 
programme COMPETE and by the Portuguese Government through the FCT --
Funda\c c\~ao para a Ci\^encia e a Tecnologia under the project
PEst-C/MAT/UI0144/2011. The first author also acknowledges the support
of the FCT project SFRH/BPD/65428/2009.

\end{document}